\documentstyle[12pt]{article}
\begin{document}
\setlength{\textheight}{574pt}
\setlength{\textwidth}{432pt}
\setlength{\oddsidemargin}{18pt}
\setlength{\topmargin}{14pt}
\setlength{\evensidemargin}{18pt}
\newtheorem{theorem}{Theorem}[section]
\newtheorem{lemma}{Lemma}[section]
\newtheorem{corollary}{Corollary}[section]
\newtheorem{conjecture}{Conjecture}
\newtheorem{remark}{Remark}[section]
\newtheorem{definition}{Definition}[section]
\newtheorem{problem}{Problem}
\newtheorem{proposition}{Proposition}[section]
\newtheorem{sublemma}{Sublemma}[section]
\title{Boundedness of {\bf Q}-Fano varieties with Picard number one}
\date{May, 1999}
\author{Hajime TSUJI}
\maketitle
\begin{abstract}
We prove birational boundedness of {\bf Q}-Fano varieties with Picard number one in arbitrary dimension.
\end{abstract}
\tableofcontents
\section{Introduction}
Let $X$ be a projective variety. 
$X$ said to be {\bf Q}-Fano variety, if 
\begin{enumerate}
\item  $X$ is {\bf Q}-factorial with only log-terminal singularities \footnote{
See Definition 5.1. In the usual definition we assume that $X$ has only terminal singularities instead of log-terminal singularities. In \cite{b}, Borisov called a {\bf Q}-Fano variety with Picard number one in this sense a Fano log-variety.},
\item $-K_{X}$ is ample.
\end{enumerate}
{\bf Q}-Fano varieties play important roles in the classification 
theory of algebraic varieties. 
About {\bf Q}-Fano varieties the following conjecture is well known.
\begin{conjecture} 
For a fixed positive integer $n$, the set of {\bf Q}-Fano $n$-folds is bounded.
\end{conjecture}
This conjecture has been solved in the case of smooth Fano varieties (\cite{n2,k-m-m} and in the case of {\bf Q}-Fano threefolds with Picard number one and 
only terminal singularities 
(\cite{ka}).   A. Borisov considered the case of {\bf Q}-Fano threefolds with Picard number one (but he also assumed a bound on the global index).

A {\bf Q}-Fano variety with Picard number one is particulary interesting 
in minimal model theory because it appears as a fiber of an 
extremal contraction (cf. \cite{ka}). 
In this paper we prove the following theorems. 

\begin{theorem}
For every positive integer $n$ the set of {\bf Q}-Fano varieties
of dimension $n$  with Picard number one is birationally bounded.  
More precisely any {\bf Q}-Fano $n$-fold with Picard number one 
is birational to a subvariety of 
degree $\leq 2^{n}n^{2n^{2}+n}$ in a projective space. 
\end{theorem}

\begin{theorem}
For every positive integer $n$ the set of Gorenstein Fano $n$-folds 
with Picard number one  
and  only log-terminal singularities is bounded. 
More precisely for any Gorenstein Fano $n$-fold $X$  with Picard number one
and log-terminal singularities  
there exists a positive integer $\nu \leq n(n+1)$ such that $\mid -\nu K_{X}\mid$ gives a 
morphism which is one to one onto its image and  
 \[
\deg \Phi_{\mid -\nu K_{X}\mid}(X) \leq 2^{n}n^{3n^{2}+n}(n+1)^{n^{2}}
\]
holds.
 
Also for any {\bf Q}-Fano $n$-fold $X$  with Picard number one
and log-terminal singularities  
there exists a positive integer $\nu \leq  n(n+1)\mbox{ind}(X)$ divisible by
$\mbox{ind}(X)$ such that 
$\mid -\nu K_{X}\mid$ gives a 
morphism which is one to one onto its image and  
 \[
\deg \Phi_{\mid -\nu K_{X}\mid}(X) \leq 2^{n}n^{2n^{2}+n}\nu^{n^{2}}
\]
holds.
\end{theorem}
The above two theorems follow from the following theorem and the techniques in 
 \cite{n2}. 
\begin{theorem}
Let $X$ be a projective variety with only log-terminal singularities. 
Assume that $-K_{X}$ is big. 
Then there exists a positive integer $\nu_{0}$ such that 
\begin{enumerate}
\item $\mid -\nu_{0}K_{X}\mid$ gives a birational rational embedding of 
$X$,
\item $X$ is covered by a family of subvarieties $\{ V\}$ 
with  $\dim V = r \geq 1$ and 
\[
\deg \Phi_{\mid -\nu_{0}K_{X}\mid}(V) \leq 2n^{r}r^{r},
\]
Moreover 
\[
\mu (V,-K_{X}) \leq 2n^{r}r^{r}
\]
holds\footnote{For the definition of $\mu$ see Definition 2.5.}. 
\end{enumerate}
\end{theorem}

The  proofs of the above theorems  do not require Mori theory or existence 
of rational curves.  
Instead of Mori theory we use the theory of multiplier ideal sheaves. 
This contrasts to the previous results \cite{n2,k-m-m}.

To illustracte our method we shall review the method in \cite{k-m-m}.
Let $X$ be a smooth Fano $n$-fold with Picard number $1$. 
By Mori thoery, there exists a dominant family of rational curves $\{ C\}$
such that 
\[
(-K_{X})\cdot C \leq n+1
\]
holds. 
Using the theory of Hilbert schemes, we see that any general two points on $X$ can be 
joined by a chain  of rational curves $C_{1},\ldots ,C_{r} (r\leq n)$ such that 
\[
(-K_{X})\cdot C_{i} \leq n+1 (1\leq i\leq r)
\]
hold. 
Next we smooth out this chain and connect these points by a single rational curve $\tilde{C}$ with 
\[
(-K_{X})\cdot \tilde{C} \leq n(n+1).
\]
By Siegel type lemma (see \cite{n2}), this implies the inequality  
\[
(-K_{X})^{n} \leq n^{n}(n+1)^{n}.
\]
Then Matsusaka type theorem (cf. \cite{a-s}) implies that $\mid -n(n+1)K_{X}\mid$ gives a one to one morphism from $X$ into a projective space such that 
\[
\deg_{\mid -n(n+1)K_{X}\mid}(X) \leq n^{2n}(n+1)^{2n} \,\,\,\,\footnote{The estimate given here is very rough and is not optimal.}.
\]
 
In the case of {\bf Q}-Fano $n$-fold $X$, the difficulty to perform the similar proof lies in  the following points :
\begin{enumerate}
\item It is not clear whether $X$ is rationally connected, 
\item Even if $X$ is rationally connected, it seems to be difficult to smooth out the chains on the singular variety, 
\item $-K_{X}$ is a ${\bf Q}$-divisor. Hence we cannot apply the Matsusaka type theorem directly.
This means that an upper bound of $(-K_{X})^{n}$ does not imply the boundedness.  
\end{enumerate}

The  essential point of our proof is to construct a dominant family of
subvarieties of   degree (birationally or regularly) bounded by a constant 
depending only on the dimension in {\bf Q}-Fano varieties by using 
multiplier ideal sheaves (Theorem 1.3\footnote{Actually we only need the special case (the case that $-K_{X}$ is ample)  of Theorem 1.3. 
In this case we do not need to use an AZD in the proof. }). 
This is similar to the existence of a dominant family of irreducible 
rational curves $\{ C\}$ with $(-K_{X})\cdot C \leq \dim X +1$, but is 
different in nature because we consider the birational rational embedding 
in Theorem 1.3 (this overcomes the third difficulty above.  
The rest of the proof is essentially same as \cite{n2}.

The method of the proofs is very similar to the one in my paper \cite{tu1} which studied
pluricanonical systems of projective varieties of general type.

The author would like to express his hearty thanks to Prof. J.P. Demailly who 
informed him the proof of Theorem 2.3.

\section{Multiplier ideal sheaves}
In this section, we shall review the basic definitions and properties
of multiplier ideal sheaves.
\subsection{Definition of multiplier ideal sheaves and Nadel's 
vanishing theorem}
\begin{definition}
Let $L$ be a line bundle on a complex manifold $M$.
A singular hermitian metric $h$ is given by
\[
h = e^{-\varphi}\cdot h_{0},
\]
where $h_{0}$ is a $C^{\infty}$-hermitian metric on $L$ and 
$\varphi\in L^{1}_{loc}(M)$ is an arbitrary function on $M$.
\end{definition}
The curvature current $\Theta_{h}$ of the singular hermitian line
bundle $(L,h)$ is defined by
\[
\Theta_{h} := \Theta_{h_{0}} + \sqrt{-1}\partial\bar{\partial}\varphi ,
\]
where $\Theta_{h_{0}} = \sqrt{-1}\bar{\partial}\partial\log h_{0}$ is 
the curvature of $h_{0}$ in the usual sense and $\partial\bar{\partial}$ is taken in the sense of a current.
The $L^{2}$-sheaf ${\cal L}^{2}(L,h)$ of the singular hermitian
line bundle $(L,h)$ is defined by
\[
{\cal L}^{2}(L,h) := \{ \sigma\in\Gamma (U,{\cal O}_{M}(L))\mid 
\, h(\sigma ,\sigma )\in L^{1}_{loc}(U)\} ,
\]
where $U$ runs open subsets of $M$.
In this case there exists an ideal sheaf ${\cal I}(h)$ such that
\[
{\cal L}^{2}(L,h) = {\cal O}_{M}(L)\otimes {\cal I}(h)
\]
holds.  We call ${\cal I}(h)$ the multiplier ideal sheaf of $(L,h)$.
If we write $h$ as 
\[
h = e^{-\varphi}\cdot h_{0},
\]
where $h_{0}$ is a $C^{\infty}$ hermitian metric on $L$ and 
$\varphi\in L^{1}_{loc}(M)$ is the weight function, we see that
\[
{\cal I}(h) = {\cal L}^{2}({\cal O}_{M},e^{-\varphi})
\]
holds.
We have the following vanishing theorem.

\begin{theorem}(Nadel's vanishing theorem \cite[p.561]{n})
Let $(L,h)$ be a singular hermitian line bundle on a compact K\"{a}hler
manifold $M$ and let $\omega$ be a K\"{a}hler form on $M$.
Suppose that $\Theta_{h}$ is strictly positive, i.e., there exists
a positive constant $\varepsilon$ such that
\[
\Theta_{h} \geq \varepsilon\omega
\]
holds.
Then ${\cal I}(h)$ is a coherent sheaf of ${\cal O}_{M}$-ideal 
and for every $q\geq 1$
\[
H^{q}(M,{\cal O}_{M}(K_{M}+L)\otimes{\cal I}(h)) = 0
\]
holds.
\end{theorem}

\subsection{Multiplier ideal sheaves on singular varieties}
Let $X$ be a normal projective variety.
Then the canonical sheaf is defined by  
\[
{\cal O}_{X}(K_{X}) := i_{*}{\cal O}_{X_{reg}}(K_{X}),
\]
where $i : X_{reg}\longrightarrow X$ is the canonical injection.
Let $L$ be a line nundle on $X$ or ${\bf Q}$-Cartier divisor such that 
$K_{X}+ L$ is Cartier (in this case $K_{X}$ must be  also {\bf Q}-Cartier). 
Let $\pi :\tilde{X}\longrightarrow X$ be a resolution of $X$. 
$h$ is said to be a singular hermitian metric on $L$, if $\pi^{*}L$ is 
a singular hermitian metric on $\pi^{*}L$. 
Clearly this definition does not depend on the choice of the resolution. 
We define the multiplier ideal sheaf ${\cal I}(h)$ of $L$ by 
\[
{\cal O}_{X}(K_{X}+L)\otimes{\cal I}(h)(U)
= \{ \eta\in \Gamma (U,{\cal O}_{X}(K_{X}+L))\mid 
h\cdot \eta\wedge\bar{\eta}\in L^{1}_{loc}(U)\} ,
\]
where $U$ runs open subsets of $X$. 
In the case of singular projective variety, to formulate vanishing theorems 
we need to use the following $L^{2}$-dualizing sheaf instead of  canonical sheaf.  
\begin{definition}
Let $X$ be a normal projective variety.  
We define the $L^{2}$-dualizing sheaf $K_{X}^{(2)}$ by
\[
K_{X}^{(2)}(U) = \{\eta \in \Gamma (U,{\cal O}_{X}(K_{X}))\mid 
\eta\wedge\bar{\eta} \in L^{1}_{loc}(U)\}.
\]
\end{definition}

\begin{proposition}
Let $X$ be a normal projective variety and let $L$ be a big {\bf Q}-line bundle. 
Let $x,x^{\prime}$ be distinct points on $X$. 
Suppose that the following conditions hold:
\begin{enumerate}
\item $K^{(2)}_{X} = K_{X}$,
\item $K_{X} + L$ is Cartier,
\item there exists a singular hermitian metric $h$  on $L$ such that 
${\cal O}_{X}/{\cal I}(h)$ has support at both  $x$ and $x^{\prime}$ and the support 
is isolated at $x$.
\item $\Theta_{h}$ dominates a positive multiple of a K\"{a}hler form induced by a projective embedding of $X$ into a projective space with the Fubini-Study K\"{a}hler form. 
\end{enumerate}
Then $H^{0}(X,{\cal O}_{X}(K_{X}+L))$ separates $x$ and $x^{\prime}$. 
\end{proposition}
{\em Proof}. 
Let $\pi : Y \longrightarrow X$ be a resolution (of singularities) such that 
the exceptional set is a divisor $D$ with normal crossings.
Let $D = \sum D_{i}$ be the irreducible decomposition. 
Let $\tau_{i}\in \Gamma (Y,{\cal O}_{Y}(D_{i}))$ be a section with divisor 
$D_{i}$. 
Let $\omega_{Y}$ be a $C^{\infty}$ K\"{a}hler form on $Y$ and we set 
\[
\omega = \omega_{Y} + c\cdot\sqrt{-1}\sum \partial\bar{\partial}\log (-\log\parallel\tau_{i}\parallel ),
\]
where $\parallel\tau_{i}\parallel$ be the hermitian norm of $\tau_{i}$ 
with respect to a $C^{\infty}$-hermitian metric on ${\cal O}_{Y}(D_{i})$ such that $\parallel\tau\parallel < 1$ and $c$ is a sufficiently small positive number. 
Then $\omega$ is a complete K\"{a}hler form on $X_{reg}$ with Poincar\'{e} growth at infinity. 

Since $\pi^{*}L$ is big, by Kodaira's lemma (\cite[Appendix]{k-o}), there exists an effective {\bf Q}-divisor $E$ such that $\pi^{*}L - E$ is ample on $Y$. 
Let $h_{E}$ be a $C^{\infty}$-hermitian metric on $\pi^{*}L - E$ and consider it as a singular hermitian metric on $L$ with strictly positive curvature.
Let us consider the singular hermitian metric 
\[
\tilde{h} := h^{1-\varepsilon}h^{\varepsilon}_{E}(\prod_{i}(-\log \parallel\tau_{i}\parallel )^{\varepsilon^{\prime}},
\]
on $L$, where $\varepsilon^{\prime} << \varepsilon$ are  sufficiently small positive numbers.
Then $\Theta_{\tilde{h}}$ is strictly positive on $Y$ and there exists a positive constant $c_{0}$ such that 
\[
\Theta_{\tilde{h}} \geq c_{0}\omega
\]
holds on $X_{reg}$.  
And if we take $\varepsilon$ (hence also $\varepsilon^{\prime}$) sufficiently small we see that ${\cal O}_{X}/{\cal I}(\tilde{h})$ has support at both $x$ and $x^{\prime}$ and the support is isolated at $x$.
Let $\sigma_{x}$ be a local generator of the locally free sheaf ${\cal O}_{X}(K_{X}+ L)$ at $x$ defined on an open neighbourhood $U$ of $x$. 
Let $\rho$ be a $C^{\infty}$-function on $X$ such that 
$\mbox{Supp}\,\rho \subset\subset U$, $0\leq \rho \leq 1$ and 
$\rho\equiv 1$ on a neighbourhood of $x$. 
We set 
\[
f := \bar{\partial}(\rho\cdot \sigma_{x}).
\]
If we take $U$ sufficiently smal, by the assumtion that $K_{X} = K_{X}^{(2)}$ holds, we see that 
\[
\int_{X}\mid f\mid^{2}\frac{\omega^{n}}{n!} < \infty  (n = \dim X)
\]
holds, where $\mid f\mid$ is the hermitian norm of $f$ with respect to $\tilde{h}$ and $\omega$.
Now we consider the equation 
\[
\bar{\partial}u = f
\]
on $X_{reg}$.
Then by the H\"{o}rmander's $L^{2}$-estimate for $\bar{\partial}$-operators on 
complete K\"{a}hler manifolds, we see that there exists a solution $u$ such that 
\[
(\sqrt{-1})^{\frac{n(n-1)}{2}}\int_{X_{reg}}\tilde{h}\cdot u\wedge\bar{u}
\leq \frac{1}{c_{0}}\int_{X}\mid f\mid^{2}\frac{\omega^{n}}{n!}
\]
holds.
Then by the $L^{2}$-condition $\sigma := \rho\cdot\sigma_{x} - u$ extends to an element of 
$H^{0}(X,{\cal O}_{X}(K_{X}+L))$ such that $\sigma (x) = \sigma_{x}(x)$ and 
$\sigma (x^{\prime}) = 0$. 
This completes the proof of Proposition 2.1. Q.E.D.  

\subsection{Analytic Zariski decomposition}
To study a big line bundle we introduce the notion of analytic Zariski
decompositions.
By using analytic Zariski decompositions, we can handle a big line bundles
like a nef and big line bundles.
\begin{definition}
Let $M$ be a compact complex manifold and let $L$ be a line bundle
on $M$.  A singular hermitian metric $h$ on $L$ is said to be 
an analytic Zariski decomposition, if the followings hold.
\begin{enumerate}
\item $\Theta_{h}$ is a closed positive current,
\item for every $m\geq 0$, the natural inclusion
\[
H^{0}(M,{\cal O}_{M}(mL)\otimes{\cal I}(h^{m}))\rightarrow
H^{0}(M,{\cal O}_{M}(mL))
\]
is an isomorphim.
\end{enumerate}
\end{definition}
\begin{remark} If an AZD exists on a line bundle $L$ on a smooth projective
variety $M$, $L$ is pseudoeffective by the condition 1 above.
\end{remark}

\begin{theorem}(\cite{tu,tu2})
 Let $L$ be a big line  bundle on a smooth projective variety
$M$.  Then $L$ has AZD. 
\end{theorem}
More generally as for the existence for general pseudoeffective line bundles, 
now we have the following theorem.

\begin{theorem}(\cite{d-p-s})
Let $X$ be a smooth projective variety and let $L$ be a pseudoeffective 
line bundle on $X$.  Then $L$ has an AZD.
\end{theorem}
{\em Proof}. Let  $h_{0}$ be a fixed $C^{\infty}$-hermitian metric on $L$.
Let $E$ be the set of singular hermitian metric on $L$ defined by
\[
E = \{ h ; h : \mbox{singular hermitian metric on $L$}, \Theta_{h} \,
\mbox{is positive}, \frac{h}{h_{0}}\geq 1 \}.
\]
Since $L$ is pseudoeffective, $E$ is nonempty (cf.\cite{d}).
We set 
\[
h_{L} = h_{0}\cdot\inf_{h\in E}\frac{h}{h_{0}},
\]
where the infimum is taken pointwise. 
Since  the supremum of a family of plurisubharmonic functions 
bounded uniformly from above is known to be again plurisubharmonic
if we modify the supremum on a set of measure $0$(i.e., it we take the uppersemicontinuous envelope) by the work of P. Lelong (\cite[p.26, Theorem 5]{l})
, we see that $h_{L}$ is also a 
singular hermitian metric on $L$ with $\Theta_{h}\geq 0$.
Suppose that there exists a nontrivial section 
$\sigma\in \Gamma (X,{\cal O}_{X}(mL))$ for some $m$ (otherwise the 
second condition in Definition 2.3 is empty).
We note that  
\[
\frac{1}{\mid\sigma\mid^{\frac{2}{m}}} 
\]
is a singular hermitian metric on $L$ with curvature 
$2\pi m^{-1}(\sigma )$, where we have considered the divisor $(\sigma )$ as a 
closed positive current. 
By the construction we see that there exists a positive constant 
$c$ such that  
\[
\frac{1}{\mid\sigma\mid^{\frac{2}{m}}} \geq c\cdot h_{L}
\]
holds. 
Hence
\[
\sigma \in H^{0}(X,{\cal O}_{X}(mL)\otimes{\cal I}(h_{L}^{m}))
\]
holds.   This means that $h_{L}$ is an AZD of $L$. 
\vspace{10mm}  Q.E.D. \\
The above construction also valid for the case of pseudoeffective line bundles on normal varieties.
Let $X$ be a normal projective variety and let $L$ be a pseudoeffective line bundle on $X$. 
Let $\pi : \tilde{X}\longrightarrow X$ be a resolution of singularities of $X$.
Then Theorem 2.3 implies that there exists an AZD $\tilde{h}_{L}$ on $\pi^{*}L$.
We set $h_{L}$ be the singular hermitian metric on $L$ induced by 
$\tilde{h}_{L}$. 
Then $h_{L}$ is an AZD of $L$\footnote{Although $h_{L}$ is defined outside of a set of 
measure $0$, this is enough.}.

\subsection{Lelong number and structure of closed positive currents}
\begin{definition}
Let $T$ be a closed positive $(1,1)$-current on  
the unit open polydisk $\Delta^{n}$ with center $O$.
Then by $\partial\bar{\partial}$-Poincar\'{e} lemma
there exists a plurisubharmonic function  $\varphi$ 
on $\Delta^{n}$ such that
\[
T = \frac{\sqrt{-1}}{\pi}\partial\bar{\partial}\varphi .
\]
We define the Lelong number $\nu (T,O)$ at $O$ by
\[
\nu (T,O) = \liminf_{x\rightarrow O}\frac{\varphi (x)}{\log \mid x\mid},
\]
where $\mid x\mid  = (\sum\mid x_{i}\mid^{2})^{1/2}$.
It is easy to see that $\nu (T,O)$ is independent of the choice of
$\varphi$ and local coordinate around $O$.
Let $V$ be a subvariety of $\Delta^{n}$. 
Then we define the Lelong number $\nu (T,V)$ by 
\[
\nu (T,V) = \inf_{x\in V}\nu (T,V). 
\]
\end{definition}
\begin{remark} More generally 
the Lelong number is defined for a closed positive
$(k,k)$-current on a complex manifold.
\end{remark}

The singular part $(\Theta_{h})_{sing}$ is analysed by the 
following theorem.

\begin{theorem}(\cite[p.53, Main Theorem]{s})
Let $T$ be a closed positive $(k,k)$-current on a complex manifold
$M$.
Then for every $c > 0$
\[
\{ x\in M\mid \nu (T,x)\geq c\}
\]
is a subvariety of codimension $\geq k$
in $M$.
\end{theorem}

The following lemma shows a rough relationship between 
the Lelong number $\nu(\Theta_{h},x)$ at $x\in X$. and the stalk of the multiplier
ideal sheaf ${\cal I}(h)_{x}$ at $x$. 

\begin{lemma}(\cite[p.284, Lemma 7]{b}\cite{b2},\cite[p.85, Lemma 5.3]{s})
Let $\varphi$ be a plurisubharmonic function on 
the open unit polydisk $\Delta^{n}$ with center $O$.
Suppose that $e^{-\varphi}$ is not locally integrable 
around $O$.
Then we have that
\[
\nu (\sqrt{-1}\partial\bar{\partial}\varphi ,O)\geq 2
\]
holds.
And if
\[
\nu (\sqrt{-1}\partial\bar{\partial}\varphi ,O) > 2n
\]
holds, then $e^{-\varphi}$ is not locally integrable.
\end{lemma}

\subsection{Volume of subvarieties}
Let $L$ be a big line bundle on a smooth projective variety
$X$. To measure the total positivity of $L$ on 
a subvariety of $X$. We define the following notion.
\begin{definition}(\cite{tu3})
Let $L$ be a big line bundle on a smooth projective variety
$X$ and let $h$ be an AZD of $L$. 
Let $Y$ be a subvariety of $X$ of dimension $r$.
We define the volume $\mu (Y,L)$ of $Y$ with respect to 
$L$ by 
\[
\mu (Y,L) := r!\limsup_{m\rightarrow\infty}m^{-r}
\dim H^{0}(Y,({\cal O}_{Y}(mL)\otimes{\cal I}(h^{m}))/tor ).
\]
\end{definition}
\begin{remark}
If we define $\mu (Y,L)$ by 
\[
\mu (Y,L) := r!\limsup_{m\rightarrow\infty}m^{-r}
\dim H^{0}(Y,{\cal O}_{Y}(mL))
\]
then it is totally different unless $Y = X$.
\end{remark}

\section{Stratification of varieties by multiplier ideal sheaves}
Now we start the proof of Theorem 1.3.  
For the proof of Theorem 1.1 and Theorem 1.2, we only need to consider the case that $-K_{X}$ is ample (hence it admits a $C^{\infty}$-hermitian metric with strictly positive curvature). 
In this case we can greatly simplify the proof below, since we do not need to consider AZD's. 

\subsection{Anticanonical rings}
Let $X$ be a  projective $n$-fold with only log-terminal singularities. 
Assume that $-K_{X}$ is big. 
For notational simplicity we denote $-K_{X}$ by $L$ in this section. 
Let $\mbox{ind}(X)$ be the global index of $X$ defined by
\[
\mbox{ind}(X) = \min \{ r\in {\bf Z}_{>0}\mid  rK_{X}\,\,\,\mbox{is Cartier.}\} .
\] 
Then by the assumption $\mbox{ind}(X)\cdot L$ is a Cartier divisor on $X$. 
For every closed point $x\in X$, there exist an affine open neighbourhood $U$ 
of $x$ and  $\mbox{ind}(X)$ sheeted abelian covering 
\[
\pi : \tilde{U} \longrightarrow U
\]
such that $K_{\tilde{U}} = \pi^{*}K_{U}$.  
This is called the canonical covering of $U$. 
We note that for every $m\geq 0$
\[
\Gamma (U,{\cal O}_{U}(-mK_{U})) \simeq \Gamma (\tilde{U},{\cal O}_{\tilde{U}}(-mK_{\tilde{U}}))^{G}
\]
holds, where $G$ denote the Galois group $Gal(\tilde{U}/U)$.
Hence we see that the direct sum of the plurianticanonical systems  
\[
R(X,-K_{X}) = \oplus_{m\geq 0}\Gamma (X,{\cal O}_{X}(-mK_{X}))
\]
is a ring under the assumption. 
We call it the anticanonical ring of $X$. 

\subsection{Construction of the stratification}

For the notational simplicity we shall assume that $X$ is smooth for the moment.In fact for the proof of Theorem 1.1, 1.3 we only need to consider the 
birationality.   
Hence if $X$ is singular, thanks to the ring structure 
of $R(X,L)$, we just need to consider on $X_{reg}$ and use Lamma 5.1 below 
to apply the $L^{2}$-vanishing theorem on $X_{reg}$ which admits a complete 
K\"{a}hler metric. 
Let $h$ be an AZD of $L$.
We set 
\[
X^{\circ} := \{ x\in X\mid  \nu (\Theta_{h}.x) = 0\} .
\]
Then we see that $X^{\circ}$ is a nonempty open subset\footnote{Actually $X^{\circ}$ contains a nonempty Zariski open subset of $X$, by Kodaira's lemma.}of $X$ in countable Zariski topology by Theorem 2.4 and  we see that for every $x\in X^{\circ}$, 
\[
{\cal I}(h^{m})_{x} = {\cal O}_{X,x}
\]
holds for every positive integer $m$ by Lemma 2.1.

Let us denote $\mu (X,L)$ by $\mu_{0}$.

\begin{lemma} Let $x,y$ be distinct points on $X^{\circ}$.  
We set 
\[
{\cal M}_{x,y} = {\cal M}_{x}\otimes{\cal M}_{y}
\]

Let $\varepsilon$ be a sufficiently small positive number.
Then 
\[
H^{0}(X,{\cal O}_{X}(mL)\otimes{\cal M}_{x,y}^{\lceil\sqrt[n]{\mu_{0}}
(1-\varepsilon )\frac{m}{\sqrt[n]{2}}\rceil})\neq 0
\]
for every sufficiently large $m$, where ${\cal M}_{x},{\cal M}_{y}$ denote the
maximal ideal sheaf of the points $x,y$ respectively.
\end{lemma}
{\em Proof of Lemma 3.1}.   
Let us consider the exact sequence:
\[
0\rightarrow H^{0}(X,{\cal O}_{X}(mL)\otimes{\cal I}(h^{m})\otimes
{\cal M}_{x,y}^{\lceil\sqrt[n]{\mu_{0}}(1-\varepsilon )\frac{m}{\sqrt[n]{2}}\rceil})
\rightarrow H^{0}(X,{\cal O}_{X}(mL)\otimes{\cal I}(h^{m}))\rightarrow
\]
\[
  H^{0}(X,{\cal O}_{X}
(mL)\otimes{\cal I}(h^{m})/{\cal M}_{x,y}^{\lceil\sqrt[n]{\mu_{0}}(1-\varepsilon )\frac{m}{\sqrt[n]{2}}\rceil}).
\]
Since 
\[
n!\limsup_{m\rightarrow\infty}m^{-n}\dim H^{0}(X,{\cal O}_{X}(mL)\otimes{\cal I}(h^{m})
/{\cal M}_{x,y}^{\lceil\sqrt[n]{\mu_{0}}(1-\varepsilon )\frac{m}{\sqrt[n]{2}}\rceil})
=
\mu_{0}(1-\varepsilon )^{n} < \mu_{0}
\]
hold, we see that Lemma 3.1 holds.  Q.E.D.

\vspace{10mm}

Let us take a sufficiently large positive integer $m_{0}$ and let $\sigma$
be a general (nonzero) element of  
$H^{0}(X,{\cal O}_{X}(m_{0}L)\otimes
{\cal M}_{x,y}^{\lceil\sqrt[n]{\mu_{0}}(1-\varepsilon )\frac{m_{0}}{\sqrt[n]{2}}\rceil})$.
We define a singular hermitian metric $h_{0}$ on $L$ by
\[
h_{0}(\tau ,\tau ) := \frac{\mid \tau\mid^{2}}{\mid \sigma\mid^{2/m_{0}}}.
\]
Then 
\[
\Theta_{h_{0}} = \frac{2\pi}{m_{0}}(\sigma )
\]
holds, where $(\sigma )$ denotes the closed positive current 
defined by the divisor $(\sigma )$.
Hence $\Theta_{h_{0}}$ is a closed positive current.
Let $\alpha$ be a positive number and let ${\cal I}_{0}(\alpha )$ denote
the multiplier ideal sheaf of $h_{0}^{\alpha}$, i.e.,
\[
{\cal I}_{0}(\alpha ) = 
{\cal L}^{2}({\cal O}_{X},(\frac{h_{0}}{h_{X}})^{\alpha})),
\]
where $h_{X}$ is an arbitrary $C^{\infty}$-hermitian metric on
$L$.
Let us define a positive number $\alpha_{0} (= \alpha_{0}(x,y))$ by
\[
\alpha_{0} := \inf\{\alpha > 0\mid ({\cal O}_{X}/{\cal I}(\alpha ))_{x}\neq 0\,\mbox{and}\, ({\cal O}_{X}/{\cal I}(\alpha))_{y}\neq 0\}.
\]
Since $(\sum_{i=1}^{n}\mid z_{i}\mid^{2})^{-n}$ is not locally integrable 
around $O\in \mbox{{\bf C}}^{n}$, by the construction of $h_{0}$, we see
that 
\[
\alpha_{0}\leq \frac{n\sqrt[n]{2}}{\sqrt[n]{\mu_{0}}(1-\varepsilon )}
\]
holds.
Then one of the following two cases occurs. \vspace{10mm} \\
{\bf Case} 1.1:  For every small positive number $\delta$, 
${\cal O}_{X}/{\cal I}(\alpha_{0}-\delta )$  has $0$-stalk 
at both $x$ and $y$. \\
{\bf Case} 1.2: For every small positive number $\delta$, 
${\cal O}_{X}/{\cal I}(\alpha_{0}-\delta )$  has nonzero-stalk 
at one of $x$ or $y$ say $y$. \vspace{10mm} \\

For the first we consider Case 1.1.
Let $\delta$ be a sufficiently small positive number and  
let $V_{1}$ be the germ of subscheme at $x$ defined by the ideal sheaf 
${\cal I}(\alpha_{0}+\delta )$.
By the coherence of ${\cal I}(\alpha ) (\alpha > 0)$, we see that 
if we take $\delta$ sufficiently small, then $V_{1}$ is independent
of $\delta$.  It is also easy to verify that $V_{1}$ is reduced if 
we take $\delta$ sufficiently small. 
In fact if we take a log resolution of 
$(X,\frac{\alpha_{0}}{m_{0}}(\sigma ))$, 
$V_{1}$ is the image of the divisor with discrepancy $-1$ 
(for example cf. \cite[p.207]{he}). 
Let $X_{1}$ be a subvariety of $X$ which defines a branch of $V_{1}$
at $x$. 
We consider the following two cases. \vspace{10mm} \\
{\bf Case} 2.1: $X_{1}$ passes through both $x$ and $y$, \\
{\bf Case} 2.2: Otherwise \vspace{10mm} \\

For the first we consider Case 2.1.
Suppose that $X_{1}$ is not isolated at $x$.  Let $n_{1}$ denote
the dimension of $X_{1}$.  Let us define the volume $\mu_{1}$ of $X_{1}$
with respect to $K_{X}$ by
\[
\mu_{1} := (n_{1})!\limsup_{m\rightarrow\infty}m^{-n_{1}}
\dim H^{0}(X_{1},({\cal O}_{X_{1}}(mL)\otimes{\cal I}(h^{m}))/tor).
\]
 
If we take $x$  sufficiently general, we may assume that $\mu_{1} > 0$ holds.
Also since $x,y\in X^{\circ}$, $X_{1}\cap X^{\circ}\neq \emptyset$ holds.

\begin{lemma} Let $\varepsilon$ be a sufficiently small positive number and let $x_{1},x_{2}$ be distinct points on $X_{1,reg}\cap X^{\circ}$. 
Then for every sufficiently large $m >1$,
\[
H^{0}(X_{1},({\cal O}_{X_{1}}(mL)\otimes{\cal I}(h^{m}))/tor \otimes
{\cal M}_{x_{1},x_{2}}^{\lceil\sqrt[n_{1}]{\mu_{1}}(1-\varepsilon )\frac{m}{\sqrt[n_{1}]{2}}\rceil})\neq 0
\]
holds. 
Moreover for any Cartier divisor $G$ on $X_{1}$,
\[
H^{0}(X_{1},({\cal O}_{X_{1}}(mL-G)\otimes{\cal I}(h^{m}))/tor \otimes
{\cal M}_{x_{1},x_{2}}^{\lceil\sqrt[n_{1}]{\mu_{1}}(1-\varepsilon )\frac{m}{\sqrt[n_{1}]{2}}\rceil})\neq 0
\]
holds for every sufficiently large $m >1$. 
\end{lemma}
{\em Proof}.  We note that ${\cal I}(h^{m})_{x_{i}} = {\cal O}_{X,x_{i}}(i = 1,2)$ hold.
Then the proof of Lemma 3.2 is identical as that of Lemma 3.1.  
\vspace{10mm} Q.E.D. \\

By Kodaira's lemma there is an effective ${\bf Q}$-divisor $E$ such
that $L- E$ is ample.
Let $\ell$ be a sufficiently large positive integer such that
\[
A := \ell (L- E)
\]
is a line bundle.
We shall replace $X^{\circ}$ by $X^{\circ} - \mbox{Supp} E$ and denote again by $X^{\circ}$.
\begin{lemma}\footnote{Here we have assumed that $X$ is smooth for simplicity. 
If $X$ has log-terminal singularities, we need to assume that $(m+1)L$ and $A$ are  Cartier as in the remark below. 
But the proof is essentially the same.
In the case that $L$ is ample, Lemma 3.3 follows from Serre's vanishing theorem.}
If we take $\ell$ sufficiently large, then 
\[
\phi_{m} : H^{0}(X,{\cal O}_{X}(mL+A)\otimes{\cal I}(h^{m}))\rightarrow 
H^{0}(X_{1},{\cal O}_{X_{1}}(mL+A)\otimes{\cal I}(h^{m}))
\]
is surjective for every $m\geq 0$.
\end{lemma}
{\bf Proof}.
Let us take a locally free resolution of the ideal sheaf ${\cal I}_{X_{1}}$
of $X_{1}$.
\[
0\leftarrow {\cal I}_{X_{1}}\leftarrow {\cal E}_{1}\leftarrow {\cal E}_{2}
\leftarrow \cdots \leftarrow {\cal E}_{k}\leftarrow 0.
\]
Then by the trivial extention of Nadel's vanishing theorem to 
the case of vector bundles,  we have :
\begin{sublemma}
If $\ell$ is sufficiently large,
\[
H^{q}(X,{\cal O}_{X}(mL+A)\otimes{\cal I}(h^{m})
\otimes{\cal E}_{j}) = 0
\]
holds for every $m\geq 1$, $q\geq 1$ and  $1\leq j\leq k$
\footnote{For more general vanishing theorem for singular hermitian vector bundles, see \cite{ca}.}.
\end{sublemma}
{\em Proof}. 
In fact if we take $\ell$ sufficiently large, we see that for every $j$, 
${\cal O}_{X}(A - K_{X})\otimes {\cal E}_{j}$ admits a $C^{\infty}$-hermitian metric $g_{j}$ such that
\[
\Theta_{g_{j}} \geq \mbox{Id}_{E_{j}}\otimes \omega
\]
holds, where $\omega$ is a K\"{a}hler form on $X$.
By   \cite[Theorem 4.1.2 and Lemma 4.2.2]{ca}, we complete the 
proof. {\bf Q.E.D.} 
\begin{remark}
Sublemma 3.1, 3.2 (below) and Lemma 3.3 also holds in the case that $X$ has log-terminal singularities.
In this case $K_{X} = K_{X}^{(2)}$ holds (cf. Lemma 5.1 below).  
Let us consider $X_{reg}$ as a complete K\"{a}hler manifold as in the proof of 
Proposition 2.1 above. 
Then if we assume that $mL-K_{X} = (m+1)L$ and $A$ are Cartier, the $L^{2}$-vanishing theorem  as in \cite[Theorem 4.1.2]{ca} on the complete K\"{a}hler manifold $X_{reg}$ immediately 
implies that 
\[
H^{q}(X,{\cal O}_{X}(mL+A)\otimes{\cal I}(h^{m})
\otimes{\cal E}_{j}) = 0
\]
holds for every such $m \geq 1$ and every $1\leq j\leq k$.
\end{remark} 
Let 
\[
p_{m} : Y_{m}\longrightarrow X
\]
be a successive blowing ups with smooth centers such that 
$p_{m}^{*}{\cal I}(h^{m})$ is locally free on $Y_{m}$\footnote{If $X$ has log-terminal singularity we take $p_{m}$ such that $Y_{m}$ is smooth.}. 
\begin{sublemma}\footnote{This sublemma also holds in the case that $X$ has 
log-terminal singularities.}
\[
R^{q}p_{m\,*}({\cal O}_{Y_{m}}(K_{Y_{m}})\otimes{\cal I}(p_{m}^{*}h^{m}))= 0
\]
holds for every $q\geq 1$ and $m\geq 1$.
\end{sublemma}
{\em Proof}. 
This sublemma follows from Theorem 2.1.
\vspace{10mm}
$\Box$  \\
By Sublemma 3.1, Sublemma 3.2 and the Leray spectral sequence,  we see that 
\[
H^{q}(Y_{m},{\cal O}_{Y_{m}}(K_{Y_{m}}+p_{m}^{*}(mL+A-K_{X}))
\otimes {\cal I}(p_{m}^{*}h^{m})
\otimes{\cal E}_{j}) = 0
\]
holds for every $q\geq 1$ and $m\geq 1$.
We note that by the definition of the multiplier ideal sheaves 
\[
p_{m\, *}({\cal O}_{Y_{m}}(K_{Y_{m}})\otimes{\cal I}(p_{m}^{*}h^{m}))= {\cal O}_{X}(K_{X})\otimes {\cal I}(h^{m})
\]
holds.  By Sublemma 3.1, Sublemma 3.2 and the Leray spectral sequence,  we see that 
\[
H^{q}(Y_{m},{\cal O}_{Y_{m}}(K_{Y_{m}}+p_{m}^{*}(mL +A - K_{X}))
\otimes {\cal I}(p_{m}^{*}h^{m})
\otimes p_{m}^{*}{\cal E}_{j}) = 0
\]
holds for every $q\geq 1$ and $m\geq 1$.
Hence 
\[
H^{1}(Y_{m},{\cal O}_{Y_{m}}(K_{Y_{m}}+ p_{m}^{*}(mL+A-K_{X})\otimes p_{m}^{*}{\cal I}(h^{m}))p_{m}^{*}{\cal I}_{X_{1}}) = 0
\]
holds. 
Hence every element of 
\[
H^{0}(Y_{m},{\cal O}_{Y_{m}}(K_{Y_{m}}+ p_{m}^{*}(mL+A-K_{X})\otimes {\cal I}(p_{m}^{*}h^{m}))\otimes {\cal O}_{Y_{m}}/p_{m}^{*}{\cal I}_{X_{1}}) 
\]
extends to an element of 
\[
H^{0}(Y_{m},{\cal O}_{Y_{m}}(K_{Y_{m}}+ p_{m}^{*}(mL+A-K_{X})\otimes {\cal I}(p_{m}^{*}h^{m}))) 
\]
Also there exists a natural map 
\[
H^{0}(X_{1},{\cal O}_{X_{1}}(mL)\otimes{\cal I}(h^{m}))
\rightarrow
\]
\[
H^{0}(Y_{m},{\cal O}_{Y_{m}}(K_{Y_{m}}+ p_{m}^{*}(mL+A-K_{X}))\otimes {\cal I}(p_{m}^{*}h^{m})\otimes {\cal O}_{Y_{m}}/p_{m}^{*}{\cal I}_{X_{1}}). 
\]
Hence we can extends every element of 
\[
p_{m}^{*}H^{0}(X_{1},{\cal O}_{X_{1}}(mL+A)\otimes{\cal I}(h^{m}))
\]
to an element of 
\[
H^{0}(Y_{m},{\cal O}_{Y_{m}}(K_{Y_{m}}+ p_{m}^{*}(mL+A-K_{X}))\otimes {\cal I}(p_{m}^{*}h^{m})) 
\]
Since 
\[
H^{0}(Y_{m},{\cal O}_{Y_{m}}(K_{Y_{m}}+ p_{m}^{*}(mL+A-K_{X})\otimes {\cal I}(p_{m}^{*}h^{m}))
\simeq 
\]
\[
 H^{0}(X,{\cal O}_{X}(mL+A)\otimes{\cal I}(h^{m}))
\]
holds by the isomorphism 
\[
p_{m\, *}({\cal O}_{Y_{m}}(K_{Y_{m}})\otimes{\cal I}(p_{m}^{*}h^{m})))= {\cal O}_{X}(K_{X})\otimes {\cal I}(h^{m}), 
\]
and Sublemma 3.2, this completes the proof of Lemma 3.3.
\vspace{10mm}{\bf Q.E.D.} \\ 

Let $\tau$ be a general section in 
$H^{0}(X,{\cal O}_{X}(A))$.

Let $m_{1}$ be a sufficiently large positive integer 
and let $\sigma_{1}^{\prime}$
be a general element of 
\[
H^{0}(X_{1},({\cal O}_{X_{1}}(m_{1}L -F)\otimes{\cal I}(h^{m_{1}}))/tor \otimes
{\cal M}_{x_{1},x_{2}}^{\lceil\sqrt[n_{1}]{\mu_{1}}(1-\varepsilon )\frac{m_{1}}
{\sqrt[n_{1}]{2}}\rceil}),
\]
where $x_{1},x_{2}\in X_{1}\cap X^{\circ}$ are distinct nonsingular points on $X_{1}$ and $F$ is an effective Cartier divisor independent of $m_{1}$ on $X_{1}$ such that every element of 
\[
H^{0}(X_{1},({\cal O}_{X_{1}}(m_{1}L+A-F)\otimes{\cal I}(h^{m_{1}}))/tor \otimes
{\cal M}_{x_{1},x_{2}}^{\lceil\sqrt[n_{1}]{\mu_{1}}(1-\varepsilon )\frac{m_{1}}
{\sqrt[n_{1}]{2}}\rceil}),
\]
extends to an element of 
\[
H^{0}(X_{1},{\cal O}_{X_{1}}(m_{1}L+A)\otimes{\cal I}(h^{m_{1}})).
\]
The existence of such $F$ can be verified as follows.
Let $f : Y\longrightarrow X$ be an embedded resolution of 
$X_{1}$.  Let $\hat{X}_{1}$ be the strict transform of $X_{1}$. 
Then an element of 
\[
H^{0}(\hat{X}_{1},{\cal O}_{\hat{X}_{1}}(f^{*}(m_{1}L+A))\otimes f^{*}{\cal I}(h^{m_{1}}))
\]
extends to an element of to
\[
H^{0}(Y,{\cal O}_{Y}(f^{*}(m_{1}L+A))\otimes f^{*}{\cal I}(h^{m_{1}})/f^{*}{\cal I}_{X_{1}}),
\]
if the element has enough zeros (which is independent of $m_{1}$) along $(f^{-1}(X_{1})-\hat{X}_{1})\cap \hat{X}_{1}$.
And also 
\[
{\cal O}_{\hat{X}_{1}}(f^{*}(m_{1}L+A))\otimes f^{*}{\cal I}(h^{m_{1}})
\]
is torsion free as a sheaf on $\hat{X}_{1}$ because $\hat{X}_{1}$ is regular and this sheaf is a subsheaf of a 
locally free sheaf. 
Hence the existence of $F$ is clear. 

By Lemma 3.2, we may assume that $\sigma_{1}^{\prime}$ is nonzero.
Then by Lemma 3.3 we see that   
\[
\sigma_{1}^{\prime}\otimes\tau\in
H^{0}(X_{1},{\cal O}_{X_{1}}(m_{1}L+A)\otimes{\cal I}(h^{m_{1}})\otimes
{\cal M}_{x_{1},x_{2}}^{\lceil\sqrt[n_{1}]{\mu_{1}}(1-\varepsilon )\frac{m_{1}}
{\sqrt[n_{1}]{2}}\rceil})
\]
extends to a section
\[
\sigma_{1}\in H^{0}(X,{\cal O}_{X}((m+\ell )L)
\otimes{\cal I}(h^{m+\ell})).
\]
Suppose that $x_{1},x_{2}\not\in \mbox{Supp}\, E$.
We may assume that  there exists a neighbourhood $U_{x_{1},x_{2}}$ of $\{ x_{1},x_{2}\}$ such that the divisor $(\sigma _{1})$  is smooth
on  $U_{x_{1},x_{2}} - X_{1}$ by Bertini's theorem, if we take $\ell$
sufficiently large because  by a trivial extention of  Lemma 3.3, 
\[
H^{0}(X,{\cal O}_{X}(mL+A)\otimes{\cal I}(h^{m}))
\rightarrow
H^{0}(X,{\cal O}_{X}(mL+A)\otimes{\cal I}(h^{m}))/
{\cal I}_{X_{1}}\cdot{\cal M}_{z})
\]
is surjective for every $z\in X$ and
 $m\geq 0$.
We define a singular hermitian metric $h_{1}$ on $L$ by
\[
h_{1} = \frac{1}{\mid\sigma_{1}\mid^{\frac{2}{m_{1}+\ell}}}.
\]
Let $\varepsilon_{0}$ be a sufficiently small positive number and 
let ${\cal I}_{1}(\alpha )$ be the multiplier ideal sheaf of 
$h_{0}^{\alpha_{0}-\varepsilon_{0}}\cdot h_{1}^{\alpha}$,i.e.,
\[
{\cal I}_{1}(\alpha ) = {\cal L}^{2}({\cal O}_{X},
h_{0}^{\alpha_{0}-\varepsilon_{0}}h_{1}^{\alpha}/
h_{X}^{(\alpha_{0}+\alpha-\varepsilon_{0})}).
\]
Suppose that $x,y$ are nonsingular points on $X_{1}\cap X^{\circ}$.
Then we set $x_{1} = x, x_{2} = y$ and define $\alpha_{1}(=\alpha_{1}(x,y))> 0$ by
\[
\alpha_{1} := \inf\{\alpha\mid ({\cal O}_{X}/{\cal I}_{1}(\alpha ))_{x}
\neq 0\,\mbox{and}\, ({\cal O}_{X}/{\cal I}_{1}(\alpha ))_{y}\neq 0\}.
\]
By Lemma 3.3 we may assume that we have taken $m_{1}$ so that  
\[
\frac{\ell}{m_{1}} \leq 
\varepsilon_{0}\frac{\sqrt[n_{1}]{\mu_{1}}}{n_{1}\sqrt[n_{1}]{2}}
\]
holds.
\begin{lemma}
\[
\alpha_{1}\leq \frac{n_{1}\sqrt[n_{1}]{2}}{\sqrt[n_{1}]{\mu_{1}}} 
+ O(\varepsilon _{0})
\]
holds.
\end{lemma}
To prove Lemma 3.4, we need the following elementary lemma.
\begin{lemma}(\cite[p.12, Lemma 6]{t})
Let $a,b$ be  positive numbers. Then
\[
\int_{0}^{1}\frac{r_{2}^{2n_{1}-1}}{(r_{1}^{2}+r_{2}^{2a})^{b}}
dr_{2}
=
r_{1}^{\frac{2n_{1}}{a}-2b}\int_{0}^{r_{1}^{-{2}{a}}}
\frac{r_{3}^{2n_{1}-1}}{(1 + r_{3}^{2a})^{b}}dr_{3}
\]
holds, where 
\[
r_{3} = r_{2}/r_{1}^{1/a}.
\]
\end{lemma}
{\em Proof of Lemma 3.4.}
Let $(z_{1},\ldots ,z_{n})$ be a local coordinate on a 
neighbourhood $U_{x}$ of $x$ in $X$ such that 
\[
U_{x} \cap X_{1} = 
\{ q\in U\mid z_{n_{1}+1}(q) =\cdots = z_{n}(q)=0\} .
\] 
We set $r_{1} = (\sum_{i=n_{1}+1}^{n}\mid z_{1}\mid^{2})^{1/2}$ and 
$r_{2} = (\sum_{i=1}^{n_{1}}\mid z_{i}\mid^{2})^{1/2}$.
Then there exists a positive constant $C$ such that 
\[
\parallel\sigma_{1}\parallel^{2}\leq 
C(r_{1}^{2}+r_{2}^{2\lceil\sqrt[n_{1}]{\mu_{1}}(1-\varepsilon )\frac{m_{1}}
{\sqrt[n_{1}]{2}}\rceil})
\]
holds on a neighbourhood of $x$, 
where $\parallel\,\,\,\,\parallel$ denotes the norm with 
respect to $h_{X}^{m_{1}+\ell}$.
We note that there exists a positive integer $M$ such that 
\[
\parallel\sigma\parallel^{-2} = O(r_{1}^{-M})
\]
holds on a neighbourhood of the generic point of $U_{x}\cap X_{1}$,
where $\parallel\,\,\,\,\parallel$ denotes the norm with respect to 
$h_{X}^{m_{0}}$. 
Then by Lemma 3.5, we have the inequality 
\[
\alpha_{1} \leq (\frac{m_{1}+\ell}{m_{1}})\frac{n_{1}\sqrt[n_{1}]{2}}{\sqrt[n_{1}]{\mu_{1}}} 
+ O(\varepsilon _{0})
\] 
holds. 
By using the fact that 
\[
\frac{\ell}{m_{1}} \leq 
\varepsilon_{0}\frac{\sqrt[n_{1}]{\mu_{1}}}{n_{1}\sqrt[n_{1}]{2}}
\]
we obtain that 
\[
\alpha_{1}\leq \frac{n_{1}\sqrt[n_{1}]{2}}{\sqrt[n_{1}]{\mu_{1}}} 
+ O(\varepsilon _{0})
\]
holds.  \vspace{10mm} Q.E.D. \\
If $x$ or $y$ is a singular point on $X_{1}$, we need the following lemma.
\begin{lemma}
Let $\varphi$ be a plurisubharmonic function on $\Delta^{n}\times{\Delta}$.
Let $\varphi_{t}(t\in\Delta )$ be the restriction of $\varphi$ on
$\Delta^{n}\times\{ t\}$.
Assume that $e^{-\varphi_{t}}$ does not belong to $L^{1}_{loc}(\Delta^{n},O)$
for every $t\in \Delta^{*}$.

Then $e^{-\varphi_{0}}$ is not locally integrable at $O\in\Delta^{n}$.
\end{lemma}
Lemma 3.6 is an immediate consequence of \cite[p. 200, Theorem]{o-t}.
Using Lemma 3.6 and Lemma 3.5, we see that Lemma 3.4 still holds
in this case by letting $x_{1}\rightarrow x$ and $x_{2}\rightarrow y$.
\vspace{10mm} \\

For the next we consider Case 1.2 and Case 2.2.  
We note that in Case 2.2 by modifying $\sigma$ a little bit 
, if necessary we may assume that
$({\cal O}_{X}/{\cal I}(\alpha_{0}-\varepsilon ))_{y}\neq 0$ 
and $({\cal O}_{X}/{\cal I}(\alpha_{0}-\varepsilon^{\prime}))_{x} = 0$ hold
for a sufficiently small positive number $\varepsilon^{\prime}$. 
For example it is sufficient to replace $\sigma$ by 
the following $\sigma^{\prime}$ constructed below.

Let $X^{\prime}_{1}$ be a subvariety which defines a branch of 
\[
\mbox{Spec}({\cal O}_{X}/{\cal I}(\alpha +\delta))
\]
at $y$.  By the assumption (changing $X_{1}$, if necessary) we may assume that $X_{1}^{\prime}$ does not 
contain $x$.  
Let $m^{\prime}$ be a sufficiently large positive integer such that 
$m^{\prime}/m_{0}$ is sufficiently small (we can take $m_{0}$ 
arbitrary large). 

Let $\tau_{y}$ be a general element of 
\[
H^{0}(X,{\cal O}_{X}(m^{\prime}L)\otimes 
{\cal I}_{X_{1}^{\prime}}), 
\]
where ${\cal I}_{X_{1}^{\prime}}$ is the ideal sheaf of 
$X_{1}^{\prime}$. 
If we take $m^{\prime}$ sufficiently large,
 $\tau_{y}$ is not identically zero. 
We set 
\[
\sigma^{\prime} = \sigma\cdot\tau_{y}. 
\] 
Then we see that the new singular hermitian metric $h_{0}^{\prime}$
defined by $\sigma^{\prime}$ satisfies the desired property.

In these cases, instead of Lemma 3.2, we use the following simpler lemma.

\begin{lemma} Let $\varepsilon$ be a sufficiently small positive number and let $x_{1}$ be a smooth point on $X_{1}\cap X^{\circ}$. 
Then for a sufficiently large $m >1$,
\[
H^{0}(X_{1},({\cal O}_{X_{1}}(mL)\otimes{\cal I}(h^{m}))/tor\otimes
{\cal M}_{x_{1}}^{\lceil\sqrt[n_{1}]{\mu_{1}}(1-\varepsilon )m\rceil})\neq 0
\]
holds.
\end{lemma}

Then taking a general $\sigma_{1}^{\prime}$ in
\[
H^{0}(X_{1},({\cal O}_{X_{1}}(m_{1}L)\otimes{\cal I}(h^{m_{1}})\otimes
{\cal M}_{x_{1}}^{\lceil\sqrt[n_{1}]{\mu_{1}}(1-\varepsilon )m_{1}
\rceil})/tor),
\]
for a sufficiently large $m_{1}$.
As in Case 1.1 and Case 2.1 we obtain a proper subvariety
$X_{2}$ in $X_{1}$ also in this case.

Let $U_{0}$ be a Zariski open subset of $X$ such that 
$\Phi_{\mid mL\mid}$ is an embedding on $U_{0}$ for some $m$.
Inductively for distinct points $x,y\in U_{0}$, we construct a strictly decreasing
sequence of subvarieties
\[
X = X_{0}(x,y)\supset X_{1}(x,y)\supset \cdots \supset X_{r}(x,y)\supset X_{r+1}(x,y) = \{ x\}\cup R_{y}\,\mbox{or}\, R_{x}\cup \{y\} ,
\]
where $R_{y}$ (resp. $R_{x}$) is a subvariety such that $x$ deos not
belong to $R_{y}$ (resp. $y$ belongs to $R_{y}$).
and invariants :
\[
\alpha_{0}(x,y) ,\alpha_{1}(x,y),\ldots ,\alpha_{r}(x,y),
\]
\[
\mu_{0},\mu_{1}(x,y),\ldots ,\mu_{r}(x,y)
\]
and
\[
n >  n_{1}> \cdots > n_{r}.
\]
By Nadel's vanishing theorem (Theorem 2.1) we have the following lemma.
\begin{lemma} 
Let $x,y$ be two distinct points on $U_{0}$. 
Then for every $m\geq \lceil\sum_{i=0}^{r}\alpha_{i}(x,y)-1\rceil$,
$\Phi_{\mid mL\mid}$ separates $x$ and $y$.
\end{lemma}
{\em Proof}. 
For simplicity let us denote $\alpha_{i}(x,y)$ by $\alpha_{i}$. 
Let us define the singular hermitian metric $h_{x,y}$ of $(m+1)L$ defined by  
\[
h_{x,y} = (\prod_{i=0}^{r-1}h_{i}^{\alpha_{i}-\varepsilon_{i}})\cdot
 h_{r}^{\alpha_{r}+\varepsilon_{r}}h^{(m-1-(\sum_{i=0}^{r-1}(\alpha_{i}-\varepsilon_{i}))- (\alpha_{r}+\varepsilon_{r})-\ell\varepsilon_{A})}h_{A}^{\varepsilon_{A}},
\]
where $h_{A}$ be a $C^{\infty}$-hermitian metric on $A$ with strictly positive 
curvature\footnote{We note that $A$ is ample. Here we have considered $h_{A}$ as a singular hermitian metric on $L$ with  strictly positive curvature. } and $\varepsilon_{A}$ is a sufficiently small positive number. 
Then we see that  ${\cal I}(h_{x,y})$ defines a subscheme of 
$X$ with isolated support around $x$ or $y$ by the definition of 
the invariants $\{\alpha_{i}\}$'s and the curvature of $h_{x,y}$ is 
strictly positive. 
Then by Nadel's vanishing theorem we see that 
\[
H^{1}(X,{\cal O}_{X}(K_{X}+(m+1)L)\otimes {\cal I}(h_{x,y})) = 0.
\]
This implies that $\Phi_{\mid mL\mid}$ separates 
$x$ and $y$.   Q.E.D. \vspace{10mm} \\

We note that for a fixed $x$, $\sum_{i=0}^{r}\alpha_{i}(x,y)$ depends on $y$.
We set
\[
\alpha (x) = \sup_{y\in U_{0}}\sum_{i=0}^{r}\alpha_{i}
\]
and let 
\[
X = X_{0}\supset X_{1}\supset X_{2}\supset\cdots X_{r} \supset
X_{r+1} = \{ x\}\cup R_{y}\,\mbox{or}\, R_{x}\cup \{ y\},
\]
be the stratification which attains $\alpha (x)$,
where $R_{y}$ (resp. $R_{x}$) is a subvariety such that $x$ deos not
belong to $R_{y}$ (resp. $y$ belongs to $R_{y}$).
In this case we call it the maximal stratification at $x$.
We see that there exists a nonempty open subset $U$ 
in countable Zariski topology 
of $X$ such that on $U$ the function $\alpha (x)$ is constant
and there exists an irreducible family of stratification 
which attains $\alpha (x)$ for every $x\in U$.

In fact this can be verified as follows.
We note that the cardinarity of
\[
\{ X_{i}(x,y)\mid  x,y \in X, x\neq y (i=0,1,\ldots )\}
\]
is uncontably many, while the cardinarity of the irreducible components
of Hilbert scheme of $X$ is countably many.
We see that for fixed $i$ and very general $x$, $\{ X_{i}(x,y)\}$ should form a family on $X$.  Similary we see that for very general $x$, we may assume that
the maximal stratification $\{ X_{i}(x)\}$ forms a family.  
 This implies the existence of $U$.   

And we may also assume that the corresponding invariants $\{\alpha_{0},
\ldots ,\alpha_{r}\}$, $\{\mu_{0},\ldots ,\mu_{r}\}$,
$\{ n = n_{0}\ldots ,n_{r}\}$ are constant on $U$.
Hereafter we denote these invariants again  by the same notations for simplicity.
The following lemma is trivial.
\begin{lemma}
\[
\alpha_{i}\leq \frac{n_{i}\sqrt[n_{i}]{2}}{\sqrt[n_{i}]{\mu_{i}}} + O(\varepsilon_{i-1})
\]
hold for $1\leq i\leq r$.
\end{lemma}

\begin{proposition}
For every 
\[
m > \lceil\sum_{i=0}^{r}\alpha_{i}-1\rceil 
\]
$\mid mL\mid$ gives a birational rational map from $X$ into 
a projective space.
\end{proposition}
Combining Proposition 3.1 and Lemma 3.9, we obtain:
\begin{corollary}
For every 
\[
m \geq \sum_{i=0}^{r}\frac{n_{i}\sqrt[n_{i}]{2}}{\sqrt[n_{i}]{\mu_{i}}}
\]
$\mid mL\mid$ gives a birational rational map from $X$ into a projective space.
\end{corollary}
\section{Proof of Theorem 1.3}
In this section we shall use the same notations and conventions as in the last section. 
Let $x,x^{\prime}$ be distinct points in $U$.
Let us denote $-K_{X}$ by $L$.
Let 
\[
X = X_{0} \supset X_{1}\supset \cdots \supset X_{r} \supset \{ x\} \cup
 R_{x^{\prime}} 
\]
be the stratification at $x,x^{\prime}$. 
Let $\mu_{i} = \mu (X_{i},L)$ and 
$n_{i} = \dim X_{i}$. 
By Lemma 3.9
\[
\alpha_{i} \leq \frac{\sqrt[n_{i}]{2}n_{i}}{\sqrt[n_{i}]{\mu_{i}}} + \varepsilon_{i-1} , \]
holds, where $\varepsilon_{i-1}$ is a positive number which can be taken 
arbitrary small.
Since $X$ is complete,  
\[
\sum_{i=0}^{r}\alpha_{i} \geq 1
\]
holds. 
In fact if 
\[
\sum_{i=0}^{r}\alpha_{i} < 1
\]
holds, there exists a nonconstant holomorphic function on 
$X$. 
This is the contradiction. 
Hence 
\[
\sum_{i=0}^{r}\frac{\sqrt[n_{i}]{2}n_{i}}{\sqrt[n_{i}]{\mu_{i}}} 
\geq 1
\]
holds. 
We set 
\[
\alpha := \max_{i} \frac{\sqrt[n_{i}]{2}n_{i}}{\sqrt[n_{i}]{\mu_{i}}}.
\]
Then 
\[
\alpha \geq \frac{1}{n}
\]
holds. 
Let $j$ be the positive integer such that  
\[
\alpha = \frac{\sqrt[n_{j}]{2}n_{j}}{\sqrt[n_{j}]{\mu_{j}}} .
\]
holds.
We set $r = n_{j}$.
Then 
\[
\mu_{j} \leq 2n^{r}r^{r}
\]
holds.  
Now we use the ring structure of $R(X,-K_{X})$. 
\begin{lemma} If $\Phi_{\mid mL\mid}\mid_{X_{j}}$ is birational rational map
onto its image, then
\[
\deg \Phi_{\mid mL\mid}(X_{j})\leq m^{n_{j}}\mu_{j}
\]
holds.
\end{lemma}
{\em Proof}.
Let $p : \tilde{X}\longrightarrow X$ be the resolution of 
the base locus of $\mid mL\mid$ and let 
\[
p^{*}\mid mL\mid = \mid P_{m}\mid + F_{m}
\]
be the decomposition into the free part $\mid P_{m}\mid$ 
and the fixed component $F_{m}$. 
Let $p_{j} : \tilde{X}_{j}\longrightarrow X_{j}$ be the resolution
of the base locus of $\Phi_{m}\mid_{X_{j}}$ 
obtained by the restriction of $p$ on $p^{-1}(X_{j})$. 
Let 
\[
p_{j}^{*}(\mid mL\mid_{X_{j}}) = \mid P_{m,j}\mid + F_{m,j}
\]
be the decomposition into the free part $\mid P_{m,j}\mid$ and 
the fixed part $F_{m,j}$.
We have
\[
\deg \Phi_{\mid mL\mid}(X_{j}) = P_{m,j}^{n_{j}}
\]
holds.
Then by the ring structure of $R(X,-K_{X})$, we have that
\[
H^{0}(X,{\cal O}_{X}(\nu P_{m}))\rightarrow 
H^{0}(X,{\cal O}_{X}(m\nu L)\otimes{\cal I}(h^{m\nu}))
\]
is injective for every $\nu\geq 1$.
Hence there exists a natural morphism
\[
H^{0}(X_{j},{\cal O}_{X_{j}}(\nu P_{m,j}))
\rightarrow 
H^{0}(X_{j},({\cal O}_{X_{j}}(m\nu L)\otimes{\cal I}(h^{m\nu}))/tor)
\]
for every $\nu\geq 1$. 
This morphism is clearly injective. 
This implies that 
\[
\mu_{j} \geq \mu (X_{j},\frac{P_{m,j}}{m})
\]
holds. 
Since $P_{m,j}$ is nef and big on $X_{j}$ we see that 
\[
\mu (X_{j},\frac{P_{m,j}}{m}) = m^{-n_{j}}P_{m,j}^{n_{j}}
\]
holds.
Hence
\[
\mu_{j}\geq m^{-n_{j}}P_{m,j}^{n_{j}}
\]
holds.  This implies that
\[
\deg \Phi_{\mid mL\mid}(X_{j})\leq \mu_{j}m^{n_{j}}
\]
holds.
\vspace{10mm} Q.E.D. \\
By Lemma 4.1, we have the following inequalities:
\[
\deg \Phi_{\mid [\alpha ]L\mid}(X_{j}) 
\leq \alpha^{r}\mu_{j} 
\leq \mu_{j}(n\frac{\sqrt[n_{j}]{2}n_{j}}{\sqrt[n_{j}]{\mu_{j}}})^{n_{j}}
= 2n^{r}r^{r}
\]
holds (we note that $\{\varepsilon_{i}\}$ can be taken arbitrary small).
Hence by the countability of the irreducible components of Hilbert scheme, 
moving $x$  as in Section 3, we complete the proof of Theorem 1.3.
\vspace{10mm} Q.E.D. \\
The following proposition is trivial by the above argument.
\begin{proposition} 
If 
\[
\mu (X,-K_{X}) \geq 2n^{2n}
\]
then $\dim V < \dim X = n$ holds. 
\end{proposition}

\section{Proofs of Theorem 1.1 and 1.2}
In this section we shall prove Theorem 1.1 and 1.2 by using 
Theorem 1.3. 
The proofs are almost indentical  as the one in \cite{n2}. 
Hence we do not repeat the full detail. 
The only difference is the use of the covering family of curves of  low degree constructed
via Theorem 1.3 instead of the covering family of  rational curves of degree less than 
$\dim X + 1$ with respect to $-K_{X}$ on a smooth Fano variety $X$ with Picard number one 
in \cite{n2}.

\subsection{Projective embedding}

\begin{definition} Let $X$ be a projective variety. 
We say that $X$ has only log-terminal singularity, if 
$X$ is {\bf Q}-Gorenstein, i.e. $K_{X}$ is a {\bf Q}-Cartier divisor  and there exists a 
resolution $\pi : \tilde{X}\longrightarrow X$ such that 
\begin{enumerate}
\item the exceptional set is a divisor $E = \sum E_{i}$ with normal crossings, 
\item $K_{\tilde{X}} = \pi^{*}K_{X} + \sum a_{i}E_{i}$ with 
$a_{i} > -1$ for all $i$.
\end{enumerate}
\end{definition}

The following lemma is an immediate consequence of the definition of log-canonical singulatities.
\begin{lemma}
Let $X$ be a projective variety with only log-terminal singularities. 
Then $K_{X} = K_{X}^{(2)}$ holds. 
\end{lemma}
The importance of Lemma 5.1 is in the fact  that we can construct only sections of 
$H^{0}_{(2)}(X_{reg},{\cal O}_{X}(mL))$ by the $L^{2}$-estimates of 
$\bar{\partial}$-operators (here we consider  sections of ${\cal O}_{X}(mL)$ 
as a $(m+1)L$-valued canonical forms on $X$) . 
Hence is is important to know whether $K_{X} = K_{X}^{(2)}$ holds or not.
\begin{lemma}
Let $X$ be a Gorenstein Fano variety with only log-terminal singularities. 
Then there exists a positive integer $\nu \leq n(n+1)$ such that 
\begin{enumerate}
\item $\mid -\nu K_{X}\mid$ gives a one to one 
morphism into a projective space,  
\item $X$ is dominated by a family of irreducible curves $\{ C\}$ 
such that 
\[
\deg \Phi_{\mid -\nu K_{X}\mid}(C) \leq 2n^{r}\nu^{r}r^{r}.
\]
for some $1\leq r \leq n$
\end{enumerate}
\end{lemma} 
{\em Proof}.
By Corollary 3.1 we see that 
for every  
\[
m \geq \sum_{i=0}^{r}\frac{\sqrt[n_{i}]{2}n_{i}}{\sqrt[n_{i}]{\mu_{i}}}
\]
$\mid mL\mid$ gives a birational rational map from $X$ into a projective space.
Since $L (= -K_{X})$ is ample Cartier divisor, we see that 
$\mu_{i} \geq 1$ holds. 
Since 
\[
\sum_{i=0}^{r}\frac{\sqrt[n_{i}]{2}n_{i}}{\sqrt[n_{i}]{\mu_{i}}} \leq n(n+1)
\]
holds, 
this implies that for every $m\geq n(n+1)$, $\mid mL\mid$ gives a birational rational map from $X$ into a projective space. 
We note that since $L$ is ample, $X^{\circ} = X_{reg}$ holds.
Hence by Proposition 2.1, the proof of Proposition 3.1 implies that 
for every $m \geq n(n+1)$, $\mid mL\mid$ gives a one to one morphism from $X_{reg}$ into a projective space \footnote{If two distinct points $x,y\in X_{reg}$ 
are special, then the stratification 
$X_{0}\supset X_{1}(x,y) \supset \cdots $ may be longer than $r+1$, but in any case the length of the stratification is less than or equal to $n+1$.}. 
Now the first assetion of Theorem 1.2 follows from the following lemma (just replace Lemma 3.6 by the lemma below) and Lemma 5.1. 
\begin{lemma} Let $Z$ be a closed $n$-dimensional subvariety of the unit open polydisk $\Delta^{N}$ and let $\varphi$ be a plusisubharmonic function on 
$Z \times \Delta$, where $\Delta$ is an open unit disk in ${\bf C}$.
Let $t$ be the standard coordinate on $\Delta$.
Then there exists a positive constant $C_{Z}$  depending only on $Z$ such that 
for every  $f\in \Gamma (Z,{\cal O}_{Z}(K_{Z}))$ such that 
\[
(\sqrt{-1})^{\frac{n(n-1)}{2}}\int_{Z}e^{-\varphi}f\wedge\bar{f} < \infty
\]
there exists a holomorphic $(n+1)$-form $F$ on 
$Z\times \Delta$ such that 
\begin{enumerate}
\item $F\mid_{Z} = dt\wedge f$,
\item $(\sqrt{-1})^{\frac{n(n+1)}{2}}\int_{Z\times\Delta}e^{-\varphi}F\wedge\bar{F}
\leq C_{Z}(\sqrt{-1})^{\frac{n(n-1)}{2}}\int_{Z}e^{-\varphi}f\wedge\bar{f}$
\end{enumerate}
\end{lemma}
This lemma is an immediate consequence of the $L^{2}$-extention theorem 
 (\cite[p. 200, Theorem]{o-t}). 

Next we shall take the dominating family of $r$-dimensional  subvarieties ${\cal V} = \{ V\}$ as in Theorem 1.3.
Let ${\cal C} := \{ C\}$ be the family of general hyperplane sections 
\[
C = H_{1}\cap \cdots \cap H_{r-1} \cap V, 
\]
where $\{ H_{i}\}$ are hyperplanes.  
Then the lemma is clear by the definition of $\{ V\}$ in Section 4 ($\{ V\}$ is a family of the strata $X_{j}$'s in the proof of Theorem 1.3) and Lemma 4.1. Q.E.D.
\subsection{Covering lemma} 
The following lemma is the key point of the argument in \cite{n2}. 
\begin{lemma}(\cite[p.687, Lemma 2]{n2})
Let $M$ be a complete irreducible variety, and let 
$N\subset M \times T$ be a family of proper subscheme of $M$ parametrized by a projective scheme $T$.
Suppose that $M$ is a union of all the $N$. 
Then there exists a closed hypersurface in $M$ that is the union of 
some of the $N_{t}$.
\end{lemma}
\subsection{Completion of the proof}
Let $\nu \leq n(n+1)$ be the positive integer as in Lemma 5.2.
Suppose that  
\[
\deg \Phi_{\mid-\nu K_{X}\mid}(X) > n^{n}(2n^{r}r^{r}\nu^{r})^{n}
\]
holds. 
We set 
\[
d := 2n^{r}r^{r}\nu^{r}.
\]
Fix any $x\in X_{reg}$. 
Then for sufficietly large $m$, there exists a section
\[
\sigma \in \Gamma (X,{\cal O}_{X}(-m\nu K_{X})) - \{ 0\}
\]
whose vanishing order at $x$ strictly exceeds $ndm$. 
We set 
\[
W_{j} = \{ p\in X\mid \mbox{mult}_{p}(\sigma ) \geq 1+mdj\} .
\]
We have a decending chain of subvarieties
\[
W_{0}\supset W_{1}\supset \cdots \supset W_{n}\ni x.
\]
For each $j\in \{ 0,\ldots ,n\}$, we choose an irreducible component of $W_{j}$ containing $x$, and denote this irreducible component by $W_{j}^{\prime}$.
We may assume that these irreducible components have been chosen so that 
\[
W_{0}^{\prime}\supset W_{1}^{\prime}\supset \cdots \supset W_{n}^{\prime}\ni x.
\]
Since this chain has length greater than $n = \dim X$, there exists some 
$j$ such that 
\[
W_{j}^{\prime} = W_{j+1}^{\prime}
\]
holds. 
We set $W = W_{j}^{\prime}$. 

Let ${\cal C}$ be the family of curves as in Lemma 5.3.
Then as in \cite[p. 687]{n2},  by the product theorem \cite[p.686, Theorem 2]{n2}(the genus of $C$ may be positive, but the product theorem can be easily generalized to the higher genus case as in \cite{c}) we see that  for every 
$C\in {\cal C}$, we see that if   
\[
C \cap W_{j+1} \neq \emptyset
\]
holds, then 
\[
C \subset W_{j}
\]
holds.
Hence by the definition of $W$, for every $C\in {\cal C}$
either
\[
C \subset W
\]
or 
\[
C \cap W = \emptyset
\]
holds.

Now we move $x$ in $X_{reg}$. 
Then we may assume that the above $W$ forms a dominating family of subvarieties
${\cal W} = \{ W_{t}\}_{t\in T}$. 
Then by the covering lemma (Lemma 5.4)we see that there exists a hypersurface $D$ which is a union of some $W_{t}$. 
Then for every $C \in {\cal C}$, we see that 
either
\[
C \subset D
\]
or 
\[
C \cap D= \emptyset
\]
holds.
Since ${\cal C}$ is a dominating family, there exists $C_{0}\in {\cal C}$ such that 
\[
C_{0}\cap D = \emptyset 
\]
holds. 
This contradicts to the assumption that the Picard number of $X$ is one. 
Hence we have that 
\[
\deg \Phi_{\mid-\nu K_{X}\mid}(X) \leq n^{n}(2n^{r}r^{r}\nu^{r})^{n}
\]
holds ($1\leq r \leq n$).
Since $1\leq r \leq n$ and $\nu \leq n(n+1)$, we see that 
\[
\deg \Phi_{\mid-\nu K_{X}\mid}(X) \leq n^{n}(2n^{2n})^{n}n^{n^{2}}(n+1)^{n^{2}}
\]
holds.
This completes the first half proof of Theorem 1.2.
The proof of the latter half of Theorem 1.2 is parallel. \vspace{10mm} Q.E.D. \\ 
The proof of Theorem 1.1 is completely parallel\footnote{In this case we take 
$\nu$ to be $\nu_{0}$.  Hence $\Phi_{\mid-\nu K_{X}\mid}$ may not be a morphism
and we need to take a modification  $\tilde{X}$ of $X$ for the resolution of the base locus. The Picard number of $\tilde{X}$ is  bigger than one in general. But this deos not matter because we just need to take the image of $D$ in $X$.} . Hence we omit it. 
 
Author's address\\
Hajime Tsuji\\
Department of Mathematics\\
Tokyo Institute of Technology\\
2-12-1 Ohokayama, Megro 152\\
Japan \\
e-mail address: tsuji@math.titech.ac.jp

\end{document}